\documentclass[12pt]{amsart}

\theoremstyle{plain}
\newtheorem{theorem}                 {Theorem}      [section]
\newtheorem{proposition}  [theorem]  {Proposition}
\newtheorem{corollary}    [theorem]  {Corollary}

\theoremstyle{definition}
\newtheorem{example}      [theorem]  {Example}
\newtheorem{remark}       [theorem]  {Remark}
\newtheorem{definition}   [theorem]  {Definition}

\def \t{\mbox{${\mathbb T}$}}
\def \r{\mbox{${\mathbb R}$}}
\def \s{\mbox{${\mathbb S}$}}

\def \h{\mbox{${\mathbb H}$}}

\def \f{\mbox{$\varphi$}}

\DeclareMathOperator{\trace}{trace}
\DeclareMathOperator{\di}{div}
\DeclareMathOperator{\grad}{grad}
\DeclareMathOperator{\Index}{index}
\DeclareMathOperator{\nul}{nullity} \DeclareMathOperator{\Div}{div}
 
\DeclareMathOperator{\riem}{Riem}
\DeclareMathOperator{\ricci}{Ricci} \DeclareMathOperator{\vol}{Vol}

\DeclareMathOperator{\isom}{Imm}
\DeclareMathOperator{\cst}{constant}

\numberwithin{equation}{section}

\begin{document}

\title{A short survey on biharmonic maps between Riemannian manifolds}

\author{S.~Montaldo}
\author{C.~Oniciuc}

\address{Universit\`a degli Studi di Cagliari\\
Dipartimento di Matematica\\
Via Ospedale 72\\
09124 Cagliari}
\email{montaldo@unica.it}

\address{Faculty of Mathematics\\ ``Al.I. Cuza'' University of Iasi\\
Bd. Carol I no. 11 \\
700506 Iasi, ROMANIA}
\email{oniciucc@uaic.ro}

\date{}



\thanks{The first author was supported by Regione Autonoma Sardegna (Italy).
The second author was supported by a CNR-NATO (Italy) fellowship,
and by the Grant At, 73/2005, CNCSIS (Romania)}

\maketitle

\section{Introduction}

Let $C^{\infty}(M,N)$ be the space of smooth maps $\phi : (M,g)\to
(N,h)$ between two Riemannian manifolds. A map  $\phi\in
C^{\infty}(M,N)$ is called {\it harmonic} if it is a critical point of the
{\it energy} functional
$$
E:C^{\infty}(M,N)\to\r, \quad E(\phi)=\frac{1}{2}\int_{M}\,
|d\phi|^2\,v_g,
$$
and is characterized by the vanishing of the first tension field
$\tau(\phi)=\trace\nabla d\phi$. In the same vein, if we denote by
$\isom(M,N)$ the space of Riemannian immersions in $(N,h)$, then a
Riemannian immersion $\phi : (M,\phi^{\ast}h)\to (N,h)$ is called
{\it minimal} if it is a critical point of the {\it volume} functional
$$
V:\isom(M,N)\to\r , \quad
V(\phi)=\frac{1}{2}\int_{M}\,v_{\phi^{\ast}h},
$$
and the corresponding Euler-Lagrange equation is $H=0$, where $H$ is
the mean curvature vector field.

If $\phi : (M,g)\to (N,h)$ is a Riemannian immersion, then it is a
critical point of the energy in $C^{\infty}(M,N)$ if and only if
it is a minimal immersion \cite{JEJHS}. Thus, in order to study
minimal immersions one can look at harmonic Riemannian immersions.

A natural generalization of harmonic maps and minimal immersions can
be given by considering the functionals obtained integrating the
square of the norm of the tension field or of the mean curvature
vector field, respectively.
More precisely:
\begin{itemize}
\item {\it biharmonic maps} are  the critical points of the {\it bienergy} functional
$$
E_2:C^{\infty}(M,N)\to\r , \quad E_2(\phi)=\frac{1}{2}\int_{M}\,
|\tau(\phi)|^2\,v_g\,;
$$
\item {\it Willmore immersions} are the critical points of the  Willmore functional
$$
W:\isom(M^2,N)\to\r , \quad W(\phi)=\int_{M^2}\, \big(|H|^2+K\big)
\,v_{\phi^{\ast}h}\,,
$$
where $K$ is the sectional curvature of $(N,h)$ restricted to the
image of $M^2$.
\end{itemize}

\noindent While the above variational problems are natural generalizations
of harmonic maps and minimal immersions, biharmonic Riemannian
immersions do not recover Willmore immersions, even when the
ambient space is $\r^n$. Therefore, the two generalizations give rise
to different variational problems.

In a different setting, in  ~\cite{BYC1}, B.Y.~Chen  defined
biharmonic submanifolds $M\subset\r^n$ of the Euclidean space as
those with harmonic mean curvature vector field, that is 
$\Delta H =0$,
where $\Delta$ is the rough Laplacian. If we apply the definition
of biharmonic maps to Riemannian immersions into the Euclidean
space we recover Chen's notion of biharmonic submanifolds. Thus
biharmonic Riemannian immersions can also be thought as a
generalization of Chen's biharmonic submanifolds.

In the last decade there has been a growing interest in the theory
of biharmonic maps which can be divided in two main research
directions. On the one side, the differential geometric aspect has
driven attention to the construction of examples and
classification results; this is the face of biharmonic maps we
shall try to report. The other side is the analytic aspect from
the point of view of PDE: biharmonic maps are solutions of a
fourth order strongly elliptic semilinear PDE. We shall not report on this
aspect and we refer the reader to ~\cite{TL1,TL2,PS,CW1,CW2} and
the references therein.

The differential geometric aspect of biharmonic submanifolds  was
also studied in the semi-Riemannian case. We shall not discuss
this case, although it is very rich in examples, and we refer the
reader to  ~\cite{BYC2} and the references therein.

We mention some other reasons that should encourage the study of
biharmonic maps.
\begin{itemize}
\item The theory of biharmonic functions is an old and rich
subject: they have been studied since 1862 by Maxwell and Airy to
describe a mathematical model of elasticity; the theory of
polyharmonic functions was later on developed, for example, by
E.~Almansi, T.~Levi-Civita and M.~Nicolescu. Recently, biharmonic
functions on Riemannian manifolds were studied by R.~Caddeo and
 L.~Vanhecke \cite{RC,RCLV},  L.~Sario
et all \cite{LS}, and others. 
\item The identity map of a
Riemannian manifold is trivially a harmonic map, but in most cases
is not stable (local minimum), for example consider $\s^{n},\,
n>2$. In contrast, the identity map, as a biharmonic map, is
always stable, in fact an absolute minimum of the energy. 
\item
Harmonic maps do not always exists, for instance, J.~Eells and
J.C.~Wood showed  in ~\cite{JEJCW} that there exists no harmonic
map from $\t^2$ to $\s^2$ (whatever the metrics chosen) in the
homotopy class of Brower degree $\pm 1$. We expect biharmonic maps
to succeed where harmonic maps have failed.
\end{itemize}

In this short survey we try to report on the theory of biharmonic
maps between Riemannian manifolds, conscious that we might have
not included all known results in the literature.

\centerline{Table of Contents}

\noindent{\bf 2.}\hspace{2mm}{\bf The biharmonic equation}\\
{\bf 3.}\hspace{2mm}{\bf Non-existence results}\\
\hspace*{6mm}{3.1.}\hspace{2mm}{Riemannian immersions}\\
\hspace*{6mm}{3.2.}\hspace{2mm}{Submanifolds of $N(c)$}\\
\hspace*{6mm}{3.3.}\hspace{2mm}{Riemannian submersions}\\
{\bf 4.}\hspace{2mm}{\bf Biharmonic Riemannian immersions}\\
\hspace*{6mm}{4.1.}\hspace{2mm}{Biharmonic curves on surfaces}\\
\hspace*{6mm}{4.2.}\hspace{2mm}{Biharmonic curves of the Heisenberg group $\h_3$}\\
\hspace*{6mm}{4.3.}\hspace{2mm}{The biharmonic submanifolds of ${\s}^3$}\\
\hspace*{6mm}{4.4.}\hspace{2mm}{Biharmonic submanifolds of $\s^n$}\\
\hspace*{6mm}{4.5.}\hspace{2mm}{Biharmonic submanifolds in Sasakian space forms}\\
{\bf 5.}\hspace{2mm}{\bf Biharmonic Riemannian submersions}\\
{\bf 6.}\hspace{2mm}{\bf Biharmonic maps between Euclidean spaces}\\
{\bf 7.}\hspace{2mm}{\bf Biharmonic maps and conformal changes}\\
\hspace*{6mm}{7.1.}\hspace{2mm}{Conformal change on the domain}\\
\hspace*{6mm}{7.2.}\hspace{2mm}{Conformal change on the codomain}\\
{\bf 8.}\hspace{2mm}{\bf Biharmonic morphisms}\\
{\bf 9.}\hspace{2mm}{\bf The second variation of biharmonic maps}
\vspace{4mm}

\noindent{\bf Acknowledgements}. The first author wishes to
thank the organizers of the ``II Workshop in Differential
Geometry  - Cordoba - June 2005'' for their exquisite
hospitality and the opportunity of presenting a lecture.
The second author wishes to thank the Dipartimento di
Matematica e Informatica,
Universit\`a di Cagliari, for hospitality
during the preparation of this paper.

\section{The biharmonic equation}\label{se-bih-eqn}
Let  $\phi : (M,g) \to (N,h)$ be a smooth map, then,
for a compact subset $\Omega\subset M$, the  {
energy} of $\phi$ is defined by
$$
E(\phi) = \frac{1}{2} \int_{\Omega}|d\phi|^{2} v_{g}=\int_{\Omega}e(\phi) v_{g}.
$$
Critical points of the energy, for any compact subset
$\Omega\subset M$, are called {harmonic maps} and the
corresponding Euler-Lagrange equation is
$$
\tau(\phi) = \trace_{g}\nabla d\phi =0.
$$
The equation $\tau(\phi)=0$ is called the {\em harmonic equation}
and, in local coordinates $\{x^i\}$ on $M$ and $\{u^{\alpha}\}$ on
$N$, takes the familiar form
$$
\tau(\phi)= \Big(-\Delta\phi^{\alpha}+
g^{ij}\,{^{N}\Gamma}^{\alpha}_{\beta\gamma}
\frac{\partial\phi^{\beta}}{\partial x^i}
\frac{\partial\phi^{\gamma}}{\partial x^j}\Big)
\frac{\partial}{\partial u^{\alpha}}=0,
$$
where ${^{N}}\Gamma^{\alpha}_{\beta\gamma}$ are the
Christoffel symbols of $(N,h)$ and  $\Delta=-\Div(\grad)$ is the
Beltrami-Laplace operator on $(M,g)$.

A smooth map $\phi : (M,g) \to (N,h)$ is {biharmonic} if it is a critical point,
for any compact subset $\Omega\subset M$, of the {bienergy} functional
$$
E_2(\phi) = \frac{1}{2} \int_{\Omega}|\tau(\phi)|^{2}\; v_{g}.
$$
We will now derive the {\em biharmonic equation}, that is the
Euler-Lagrange equation associated to the bienergy. For simplicity
of exposition we will perform the calculation for smooth maps $\phi: (M,g) \to \r^n$, defined by
$\phi(p)=(\phi^1(p),\ldots,\phi^n(p))$, with $M$ compact. In this
case we have
\begin{equation}\label{eq:beltrami}
\tau(\phi)=-\Delta \phi = -(\Delta\phi^1,\ldots,\Delta\phi^n)\quad\text{and}\quad
E_2(\phi)=\frac{1}{2} \int_M |\Delta \phi|^2 v_g\,.
\end{equation}
To compute the corresponding Euler-Lagrange equation, let
$\phi_t=\phi+t X$ be a one-parameter smooth variation of $\phi$ in
the direction of a vector field $X$ on $\r^n$ and denote with
$\delta$ the operator ${d}/{dt}|_{t=0}$. We have
\begin{align*}
\delta (E_2(\phi_t)) &=
\,\int_{M^2} \langle\delta \Delta \phi_t,\Delta \phi\rangle v_g
=\,\int_{M^2} \langle\Delta X, \Delta \phi\rangle v_g \\
&=\,\int_{M^2} \langle X,\Delta^2\phi\rangle v_g\,,
\end{align*}
where in the last equality we have used that  $\Delta$ is
self-adjoint. Since $\delta (E_2(\phi_t))=0$, for any vector field
$X$,  we conclude that $\phi$ is biharmonic if and only if
$$
\Delta ^2\phi=0.
$$
Moreover, if $\phi:M\to\r^n$ is a Riemannian immersion, then, using Beltrami equation $\Delta \phi= - m H$,  we have that $\phi$ is biharmonic if and only if
$$
\Delta ^2\phi=- m \Delta H=0.
$$
Therefore, as mentioned in the introduction, we recover Chen's
definition of biharmonic submanifolds in $\r^n$.

For a smooth map $\phi : (M,g) \to (N,h)$ the Euler-Lagrange
equation associated to the bienergy becomes more complicated and,
as one would expect, it involves the curvature of the codomain.
More precisely, a smooth map $\phi : (M,g) \to (N,h)$ is
biharmonic if it satisfies the following biharmonic equation
$$
\tau_{2}(\phi) = - \Delta^{\phi} \tau(\phi) -
\trace_{g}{R^{N}}(d\phi,\tau(\phi))d\phi =0,
$$
where $\Delta^{\phi} = - \trace_{g} \big(
\nabla^{\phi}\nabla^{\phi} - \nabla^{\phi}_{\nabla} \big)$ is the
rough Laplacian on sections of $\phi^{-1}TN$ and
$R^{N}(X,Y)=[\nabla_X,\nabla_Y]-\nabla_{[X,Y]}$ is the
curvature operator on $N$.

\noindent From the expression of the {\em bitension} field
$\tau_2$ it is clear that a harmonic map ($\tau=0$) is
automatically a biharmonic map, in fact a minimum of the bienergy.

We call a non-harmonic biharmonic map a {\it proper} biharmonic
map.

\section{Non-existence results}
As we have just seen, a harmonic map is biharmonic, so a basic
question in the theory is to understand under what conditions the
converse is true. A first general answer to this problem,
proved by G.Y.~Jiang, is
\begin{theorem}[\cite{GYJ1,GYJ2}]\label{jian-nonex}
Let $\phi:(M,g)\to (N,h)$ be a smooth map. If $M$
is compact, orientable and $\riem^N\leq 0$, then $\phi$ is
biharmonic if and only if it is harmonic.
\end{theorem}
Jiang's theorem is a direct application of the Weitzenb\"ock formula.
In fact, if $\phi$ is biharmonic, the  Weitzenb\"ock formula and
$\tau_2(\phi)=0$ give
\begin{eqnarray*}
\frac{1}{2}\Delta\vert\tau(\phi)\vert^2&=&\langle\Delta\tau(\phi),
\tau(\phi)\rangle-\vert d\tau(\phi)\vert^2 \\
&=&\trace\langle
R^N(\tau(\phi),d\phi)d\phi,\tau(\phi)\rangle-\vert
d\tau(\phi)\vert^2 \leq 0.
\end{eqnarray*}
Then, since $M$ is compact, by the maximal principle, we find that
$d\tau(\phi)=0$. Now using the identity
$$
\di\langle d\phi,\tau\rangle=\vert\tau(\phi)\vert^2+\langle d\phi,d\tau(\phi)\rangle,
$$
we deduce that
$\di\langle d\phi,\tau\rangle=\vert\tau(\phi)\vert^2$ and, after
integration, we conclude.

\subsection{Riemannian immersions}
If $M$ is not compact, then the above argument can be used with
the extra assumption that $\phi$ is a Riemannian immersion and
that the norm of $\tau(\phi)$ is constant, as was shown
by C.~Oniciuc in
\begin{theorem}[\cite{CO1}]\label{oniciuc-non-ex}
 Let $\phi:(M,g)\to
(N,h)$ be a Riemannian immersion. If $\vert\tau(\phi)\vert$ is
constant and $\riem^N\leq 0$, then $\phi$ is biharmonic if and
only if it is minimal.
\end{theorem}

The curvature condition in Theorem~\ref{jian-nonex}
and~\ref{oniciuc-non-ex} can be weakened in the case of
codimension one, that is $m=n-1$. We have
\begin{theorem}[\cite{CO1}]
Let $\phi:(M,g)\to (N,h)$ be a Riemannian immersion with $\ricci^N\leq 0$ and $m=n-1$.
\begin{itemize}
\item[a)] If $M$ is compact and orientable, then $\phi$ is biharmonic
if and only if it is minimal.
\item[b)]  If $\vert\tau(\phi)\vert$ is
constant, then $\phi$ is biharmonic if and only if it is minimal.
\end{itemize}
\end{theorem}

\subsection{Submanifolds of $N(c)$}

\noindent Let $N(c)$ be a manifold with constant sectional curvature
$c$, $M$  a submanifold of $N(c)$ and denote by ${\bf
i}:M\to N(c)$ the canonical inclusion.
In this case the tension and bitension fields assume the
following form
$$
\tau({\bf i})=mH, \quad \tau_2({\bf i})=-m\big(\Delta H - m c H\big).
$$

If $c\leq 0$, there are strong restrictions on the existence of proper
 biharmonic submanifolds in $N(c)$.
If $M$ is compact, then there exists no proper biharmonic Riemannian
immersion from $M$ into $N(c)$. In fact, from
Theorem~\ref{jian-nonex}, $M$ should
 be minimal.
If $M$ is not compact and ${\bf i}$ is a proper biharmonic map
then, from Theorem~\ref{oniciuc-non-ex}, $|H|$ cannot be constant.

If $c> 0$, as we shall see in Section~\ref{sec-s3}
and~\ref{section-sn}, we do have examples of compact proper
biharmonic submanifolds.

The main tool in the study of biharmonic submanifolds of $N(c)$ is the 
decomposition
of the bitension field in its tangential and normal components.
 Then, asking that
both components are identically zero, we conclude that  the canonical
inclusion ${\bf i}:M\to N(c)$ is biharmonic if and
only if
\begin{equation}
\label{eq:biharmonic-c-curvature}
\left\{
\begin{array}{cll}
\Delta^{\perp}H+\trace B(\cdot,A_H\cdot)- c\; m H&=&0 \\
4\trace A_{\nabla^{\perp}_{(\cdot)}H}(\cdot)+m\grad (\vert H
\vert^2)&=&0,
\end{array}
\right.
\end{equation}
where $B$ is the second fundamental form of $M$ in $N(c)$, $A$ the
shape operator, $\nabla^{\perp}$ the normal connection and
$\Delta^{\perp}$ the Laplacian in the normal bundle of $M$.

Equation~\eqref{eq:biharmonic-c-curvature} was used by B.Y.~Chen,
for $c=0$, and by R.~Caddeo, S.~Montaldo and C.~Oniciuc, for
$c<0$, to prove that in the case of biharmonic surfaces in
$N^3(c),\;c\leq 0$, the mean curvature must be constant, thus
\begin{theorem}
[\cite{RCSMCO2},\cite{BYCIS}] \label{eq:biharmonic-in-r3}Let $M^2$ be
a surface of $N^3(c),\;c\leq 0$. Then $M$ is biharmonic if and only if it is minimal.
\end{theorem}

For higher dimensional cases it is not known whether there
exist proper biharmonic submanifolds of $N^n(c),\;n>3,\;c\leq 0$,
although, for $N^n(c)=\r^n$, partial results have been obtained.
For instance:
\begin{itemize}
\item  Every biharmonic curve of $\r^n$ is an open part of a
straight line \cite{ID}.
\item Every biharmonic submanifold of finite type in $\r^n$ is
minimal \cite{ID}.
\item There exists no proper biharmonic hypersurface of $\r^n$ with at
most two principal curvatures \cite{ID}.
\item Let $M^m$ be a pseudo-umbilical submanifold of $\r^n$.
If $m\neq 4$, then $M$ is biharmonic if and only if  minimal \cite{ID}.
\item Let $M^3$ be a hypersurface of $\r^4$. Then $M$ is biharmonic if and only if
 minimal \cite{THTV}.
\item Let $M$ be a submanifold of $\s^n$. Then it is biharmonic in $\r^{n+1}$
if and only if  minimal in $\r^{n+1}$ \cite{BYC1}.
\item Let $M^m$ be a pseudo-umbilical submanifold of $N(-1)$.
If $m\neq 4$, then $M$ is biharmonic if and only if minimal \cite{RCSMCO2}.
\end{itemize}

All this results suggested the following

\noindent {\bf Generalized Chen's Conjecture}:
{\it Biharmonic submanifolds of a manifold $N$ with $Riem^{N}\leq 0$ are
minimal}.

\subsection{Riemannian submersions}

Let $\phi:(M,g)\to (N,h)$ be a Riemannian submersion with basic
tension field. Then the bitension field, computed in \cite{CO1}, is
\begin{equation}\label{eq:bihar-submersion}
\tau_2(\phi)=\trace
^N\nabla^2\tau(\phi)+^N\nabla_{\tau(\phi)}\tau(\phi)+
\ricci^N\tau(\phi).
\end{equation}
Using this formula we find some non-existence results
which are, in some sense, dual to those for Riemannian immersions. They can
be stated as follows:
\begin{proposition}[\cite{CO1}]
A biharmonic Riemannian submersion $\phi:M\to N$ with basic tension field
is harmonic in the following cases:
\begin{itemize}
\item[a)] if $M$ is compact, orientable and $\ricci^N\leq 0$;
\item[b)] if $\ricci^N<0$ and $\vert\tau(\phi)\vert$ is constant;
\item[c)] if $N$ is compact and $\ricci^N<0$.
\end{itemize}
\end{proposition}

\section{Biharmonic Riemannian immersions}
In this section we report on the known examples of proper 
biharmonic Riemannian immersions.
Of course, the first and easiest examples can be found looking at differentiable
curves in a Riemannian manifold.
This is the first class we shall describe.

Let  $\gamma:I\to (N,h)$ be a curve  parametrized by arc length from an open
interval $I\subset\r$ to a Riemannian manifold. In this case the tension field becomes
$ \tau({\gamma})=\nabla_{T}T,\,T=\gamma'$, and the
biharmonic equation reduces to
\begin{equation}\label{bieq-curve}
\nabla^3_{T}T - R(T,\nabla_T T)T = 0.
\end{equation}
To describe geometrically Equation \eqref{bieq-curve} let recall
 the definition of the Frenet frame.
\begin{definition}[See, for example, \cite{DL}] \label{def2.1}
The Frenet frame $\{F_{i}\}_{i=1,\dots,n}$ associated to a curve
$\gamma : I\subset {\r}\to (N^{n},h)$, parametrized by arc length,
 is the orthonormalisation of
the $(n+1)$-uple  $\{ \nabla_{\frac{\partial}{\partial
t}}^{(k)} d\gamma
(\frac{\partial}{\partial t})
\}_{k=0,\dots,n}$, described by:
$$
\left\{
\begin{array} {lcl}
F_{1}&=&d\gamma
(\frac{\partial}{\partial t}) ,  \\
\nabla_{\frac{\partial}{\partial
t}}^{\gamma} F_{1} &=& k_{1} F_{2} , \\
\nabla_{\frac{\partial}{\partial
t}}^{\gamma} F_{i} &=& - k_{i-1} F_{i-1}
+ k_{i}F_{i+1} , \quad \forall i =
2,\dots,n-1 , \\
\nabla_{\frac{\partial}{\partial
t}}^{\gamma} F_{n} &=& - k_{n-1} F_{n-1}
\end{array}
\right.
$$
where the functions $\{k_{1}=k>0,k_{2}=-\tau,k_{3},\ldots,k_{n-1}\}$
 are called the curvatures of $\gamma$ and $\nabla^{\gamma}$ is the connection
 on the pull-back bundle $\gamma^{-1}(TN)$. Note that
$F_{1}=T={\gamma'}$ is the unit tangent vector field along the curve.
\end{definition}
We point out  that when the dimension of $N$ is $2$, the first curvature $k_1$ is
replaced by the signed curvature.

\noindent Using the Frenet frame, we get that a curve is proper ($k_{1}\neq0$) biharmonic if and only if

\begin{equation}\label{euler-lagrange}
\left\{
\begin{array}{l}
k_{1} = \cst\neq 0\\
k_{1}^2 + k_2^2 =  R(F_1,F_2,F_1,F_2) \\
{k'_2} = - R(F_1,F_2,F_1,F_3) \\
k_2 k_3 = - R(F_1,F_2,F_1,F_4) \\
R(F_1,F_2,F_1,F_j) = 0 \hspace{1,5 cm} j = 5,\ldots, n
\end{array}
\right.
\end{equation}

\subsection{Biharmonic curves on surfaces}

Let $(N^2,h)$ be an oriented surface and let $\gamma:I\to (N^2,h)$ be
a differentiable curve parametrized by arc length. Then
Equation~\eqref{euler-lagrange} reduces to
$$
\left\{
\begin{array}{l}
k_{g} = \cst\neq 0\\
k_{g}^2=G
\end{array}
\right.
$$
where $k_{g}$ is the curvature (with sign)
 of $\gamma$ and $G=R(T,N,T,N)$ is the Gauss curvature of the surface.

As an immediate consequence we have:
\begin{proposition}[\cite{RCSMPP1}]
\label{pro-basic} Let $\gamma:I\to (N^2,h)$ be a proper
biharmonic curve on an oriented surface $N^2$. Then, along $\gamma$, the Gauss
curvature must be constant, positive and equal to the square of the
geodesic curvature of $\gamma$. Therefore, if $N^2$ has non-positive
Gauss curvature, any biharmonic curve is a geodesic of $N^2$.
\end{proposition}

\noindent Proposition~\ref{pro-basic} gives a positive answer to the generalized Chen's conjecture.

Now, let $\alpha(u)=(f(u),0,g(u))$ be a curve in the $xz$-plane and
consider the surface of revolution, obtained by rotating this curve
about the $z$-axis, with the standard parametrization
$$
X(u,v)=(f(u)\cos(v),f(u)\sin(v),g(u))\,,
$$
where $v$ is the rotation angle. Assuming that
$\alpha$ is parametrized by arc length, we have
\begin{proposition}[\cite{RCSMPP1}]
A parallel $u=u_0=\cst$ is biharmonic if and only if $u_0$
satisfies the equation
$$
{f'}^2(u_0)+{f''}(u_0)f(u_0)=0.
$$
\end{proposition}

\begin{example}[Torus]
On a torus of revolution  with its standard parametrization
$$
X(u,v)=\Big(\big(a+r\cos(\tfrac{u}{r})\big)\cos
v,\big(a+r\cos(\tfrac{u}{r})\big)\sin v, r\sin(\tfrac{u}{r}) \Big), \quad a>r,
$$
the biharmonic parallels are
$$
u_1=r\arccos\Big(\frac{-a+\sqrt{a^2+8r^2}}{4r}\Big),\quad
u_2=2r\pi-r\arccos\Big(\frac{-a+\sqrt{a^2+8r^2}}{4r}\Big).
$$
\end{example}

\begin{example}[Sphere]
There is a geometric way to understand the behaviour of
biharmonic curves on a sphere. In fact, the torsion  $\tau$ and curvature
$k$ (without sign) of
$\gamma$, seen in the ambient space $\r^3$, satisfy
$k_g({k'}_g+\tau k^2 r) = 0$. From this we see
that $\gamma$ is a proper biharmonic curve if and only if $\tau=0$ and
$k={\sqrt{2}}/{r}$, i.e. $\gamma$ is the circle of radius
${r}/{\sqrt{2}}$.
\end{example}

\noindent For more examples see \cite{RCSMPP1,RCSMPP2}.

\subsection{Biharmonic curves of the Heisenberg group $\h_3$}

The Heisenberg group $\h_3$ can be seen as the Euclidean space
$\r^3$ endowed with the multiplication
$$
(\widetilde{x},\widetilde{y},\widetilde{z})(x,y,z)=
(\widetilde{x}+x,\widetilde{y}+y,\widetilde{z}+z+
\frac{1}{2}\widetilde{x}y- \frac{1}{2}\widetilde{y}x)
$$
and with the left-invariant Riemannian metric $g$ given by
\begin{equation}\label{metricaheis}
g=dx^2+dy^2+(dz+\frac{y}{2}dx-\frac{x}{2}dy)^2.
\end{equation}
Let $\gamma:I\to \h_3$ be a differentiable curve
parametrized by arc length. Then, from \eqref{euler-lagrange}, $\gamma$
is a proper biharmonic curve if and only if
\begin{equation}\label{3.3}
\left\{
\begin{array}{l}
k=\cst\neq 0 \\
k^2+\tau^2=\frac{1}{4}-B_3^2 \\
\tau'=N_3B_3,
\end{array}
\right.
\end{equation}
where $T = T_1e_1 + T_2e_2 + T_3e_3$, $N = N_1e_1 + N_2e_2 +
N_3e_3$, and $B = T \times N = B_1e_1 + B_2e_2 + B_3e_3$.
Here $\{e_1,e_2,e_3\}$ is the left-invariant orthonormal basis with respect to
the metric \eqref{metricaheis}.

By analogy with  curves in $\r^3$, we use
the name {\it helix} for a curve in a Riemannian manifold having
 both geodesic curvature and geodesic torsion constant.

\noindent Using System~\eqref{3.3}, in \cite{RCCOPP},
R.~Caddeo, C.~Oniciuc and P.~Piu showed that a proper biharmonic curve in $\h_3$
is a helix and give their explicit parametrizations,
as shown in the following
\begin{theorem}[\cite{RCCOPP}]
\label{eq:biharmonicintegralcurves} The parametric equations of all
proper biharmonic curves $\gamma$ of $\h_3$ are
\begin{equation}
\label{eq:biharmonichelices1}
\left\{
\begin{array}{lll}
x(t) &=& \frac{1}{A}\sin\alpha_0\sin(At+a) + b, \\
y(t) &=& -\frac{1}{A}\sin\alpha_0\cos(At+a) + c, \\
z(t) &=& (\cos\alpha_0+\frac{(\sin\alpha_0)^2}{2A})t  \\
  && -\frac{b}{2A}\sin\alpha_0\cos(At + a)-
\frac{c}{2A}\sin\alpha_0\sin(At + a) + d,
\end{array}
\right.
\end{equation}
where $2A =\cos\alpha_0\pm\sqrt{5(\cos\alpha_0)^2-4}$,
$\alpha_0\in(0,\arccos\frac{2\sqrt{5}}{5}]\cup[\arccos(-\frac{2\sqrt{5}}{5})
,\pi)$ and $a,b,c,d \in \r$.
\end{theorem}

\noindent Geometrically, proper biharmonic curves in $\h_3$ can be
obtained as the intersection of a minimal helicoid with a round cylinder.
Moreover, they are geodesic of this round cylinder.

The above method can be extended to study
biharmonic curves in Cartan-Vranceanu three-manifolds ($N^3, ds^{2}_{m,\ell}$), where $N=\r^3$ if $m\geq 0$,
 $N=\{(x,y,z)\in\r^3\;:\; x^{2} + y^{2} < - \frac{1}{m}\}$ if $m<0$, and
 the Riemannian metric  $ds^{2}_{m,\ell}$ is defined by
\begin{equation}
\label{5.35}
\hspace{0.7 cm} ds^{2}_{m,\ell} =\frac{dx^{2} + dy^{2}}{[1
+ m(x^{2} + y^{2})]^{2}} +  \left(dz + \frac{\ell}{2} \frac{ydx -
xdy}{[1 + m(x^{2} + y^{2})]}\right)^{2}, \quad \ell,m \in \r .
\end{equation}
This two-parameter family of metrics reduces to the Heisenberg
metric for $m = 0 $ and $\ell = 1$.
The system for proper biharmonic curves corresponding to
the metric $ds^2_{m,\ell}$ can be obtained by using the same
techniques, and turns out to be
\begin{equation}
\label{3.3Vran} \left\{
\begin{array}{l}
k=\cst\neq 0 \\
k^2+\tau^2=\frac{\ell^2}{4}-(\ell^2 - 4m) B_3^2 \\
\tau'= (\ell^2 - 4m)N_3B_3 .
\end{array}
\right.
\end{equation}

\noindent System~\ref{3.3Vran} also implies that the proper biharmonic
curves of $(N,ds^2_{m,\ell})$ are helices \cite{RCSMCOPP}.
The explicit parametrization of proper biharmonic curves of $(N,ds^2_{m,\ell})$
was given in~\cite{JTCJIJEL}, for $\ell=1$, and in ~\cite{RCSMCOPP} in  general.

We point out that biharmonic curves were studied in other spaces which
are generalizations of the above cases. For example:
\begin{itemize}
\item In ~\cite{DF}, D.~Fetcu studied biharmonic curves in the
$(2n+1)-$dimensional Heisenberg group $\h_{2n+1}$ and obtained two
families of proper biharmonic curves.
\item A.~Balmu\c s studied, in ~\cite{AB2}, the biharmonic curves on
Berger spheres $\s^3_{\varepsilon}$, obtaining their explicit
parametric equations.
\end{itemize}

\subsection{The biharmonic submanifolds of $\s^3$}\label{sec-s3}

In \cite{RCSMCO1} the authors give a complete classification
of the proper biharmonic submanifolds of $\s^3$.

Using System\eqref{euler-lagrange} it was first proved
that the proper biharmonic curves $\gamma:I\to \s^3$ are the helices
with $k^2+\tau^2=1$. If we look at $\gamma$ as a curve
in $\r^4$, the biharmonic condition can be expressed as
\begin{equation}\label{eq-bicur}
\gamma^{\rm \i v}+2\gamma''+(1-k^2)\gamma=0.
\end{equation}

\noindent Now, by integration of  \eqref{eq-bicur}, we obtain

\begin{theorem}[\cite{RCSMCO1},\cite{ABCO}]
Let $\gamma:I\to \s^3$ be a curve parametrized by arc length. Then
it is proper biharmonic if and only if it is either the circle
of radius $\frac{1}{\sqrt{2}}$, or a geodesic of the Clifford torus
$\s^1(\frac{1}{\sqrt{2}})\times\s^1(\frac{1}{\sqrt{2}})\subset\s^3$
with slope different from $\pm 1$.
\end{theorem}

As to proper biharmonic surfaces $M^2\subset\s^3$ of the three-dimensional
sphere, one can first prove that Equation~\eqref{eq:biharmonic-c-curvature}
implies the following
\begin{theorem}[\cite{RCSMCO1}]
Let $M$ be a surface of $\s^3$. Then it is proper biharmonic  if
and only if $|H|$ is constant and $|B|^2=2$.
\end{theorem}

The classification of constant mean curvature surfaces in $\s^3$ with $|B|^2=2$
is known, in fact we have

\begin{theorem}[\cite{RCSMCO1},\cite{ZHH}]
Let $M$ be a surface of $\s^3$ with constant mean curvature and
$|B|^2=2$.
\begin{itemize}
\item[a)] If $M$ is not compact, then locally it is a piece of
either a hypersphere $\s^2(\frac{1}{\sqrt{2}})$ or a torus
$\s^1(\frac{1}{\sqrt{2}})\times \s^1(\frac{1}{\sqrt{2}})$.
\item[b)] If $M$ is compact and orientable, then it is either
$\s^2(\frac{1}{\sqrt{2}})$ or $\s^1(\frac{1}{\sqrt{2}})\times
\s^1(\frac{1}{\sqrt{2}})$.
\end{itemize}
\end{theorem}

\noindent Now, since the Clifford torus
$\s^1(\frac{1}{\sqrt{2}})\times \s^1(\frac{1}{\sqrt{2}})$ is minimal
in $\s^3$, we can state:

\begin{theorem}[\cite{RCSMCO1}]\label{class-bisurf-s3}
Let $M$ be a proper biharmonic surface of $\s^3$.
\begin{itemize}
\item[a)] If $M$ is not compact, then it is locally a piece of
$\s^2(\frac{1}{\sqrt{2}})\subset\s^3$.
\item[b)] If $M$ is compact and orientable, then it is $\s^2(\frac{1}{\sqrt{2}})$.
\end{itemize}
\end{theorem}

\subsection{Biharmonic submanifolds of $\s^n$}\label{section-sn}

We start describing some basic examples of proper biharmonic submanifolds of
 $\s^n$.

Let $\phi_t:\s^m\to \s^{m+1}$, $\phi_t(x)=(tx,\sqrt{1-t^2})$, $t\in
[0,1]$. Up to a homothetic transformation, $\phi_t$ is the canonical
inclusion of the hypersphere $\s^m(t)$ in $\s^{m+1}$.
A simple calculation shows that $E_2(\phi_t)=\frac{m^2}{2}t^2(1-t^2)\vol(\s^m)$.
Derivating $E_2(\phi_t)$ with respect to $t$ we find that
$\big(E_2(\phi_t)\big)'=0$ if and only if $t={1}/{\sqrt{2}}$.

This simple argument shows that $\s^m(a)$ is a good candidate for proper
 biharmonic submanifold of $\s^{m+1}$ if $a={1}/{\sqrt{2}}$. It is not
 difficult to show that, indeed,  the bitension field
of $\s^m({1}/{\sqrt{2}})$ is zero, proving that it is the only proper biharmonic
hypersphere of $\s^{m+1}$.

To explain the next example we first note that,
from \eqref{eq:biharmonic-c-curvature}, we have

\begin{proposition}
Let $M^{m}$ be a non-minimal hypersurface of $\s^{m+1}$ with parallel mean curvature,
i.e. the norm of $H$ is constant. Then $M^{m}$ is a
proper biharmonic submanifold if and only if
$|B|^2=m$.
\end{proposition}

\noindent Let $m_1,m_2$ be two positive integers such that $m=m_1+m_2$, and
let $r_1,r_2$ be two positive real numbers such that
$r_1^2+r_2^2=1$. Then the generalized Clifford torus
$\s^{m_1}(r_1)\times \s^{m_2}(r_2)$ is a hypersurface of $\s^{m+1}$.
A simple calculation shows that
$$
|H|=\frac{1}{m\, r_1 r_2}|m_2\, r_1^2 - m_1\, r_2^2|\quad {\rm and}
\quad|B|^2=m_1\big(\tfrac{r_2}{r_1}\big)^2+m_2\big(\tfrac{r_1}{r_2}\big)^2.
$$
We thus have
\begin{example}[\cite{GYJ1,GYJ2}]\quad
\begin{enumerate}
\item If $m_1\neq m_2$, then $\s^{m_1}(r_1)\times \s^{m_2}(r_2)$ is a
proper biharmonic submanifold of $\s^{m+1}$ if and only if
$r_1=r_2=\frac{1}{\sqrt{2}}$.
\item If $m_1=m_2=q$, then the following statements are equivalent:
\begin{itemize}
\item $\s^q(r_1)\times \s^q(r_2)$ is a biharmonic submanifold of $\s^{2q+1}$
\item $\s^q(r_1)\times \s^q(r_2)$ is a minimal submanifold of $\s^{2q+1}$
\item $r_1=r_2=\frac{1}{\sqrt{2}}$.
\end{itemize}
\end{enumerate}
\end{example}

The submanifolds $\s^m(\frac{1}{\sqrt{2}})$ and the generalized Clifford
torus are the only known examples of proper biharmonic hypersurfaces
of $\s^{m+1}$. As we have seen in Theorem~\ref{class-bisurf-s3},
for $\s^3$, the hypersphere $\s^2(\frac{1}{\sqrt{2}})$ is the only one.

{\bf Open problem:} {\it classify all proper biharmonic hypersurfaces of $\s^{m+1}$}.

The situation seems much richer if the codimension is greater than one.
We shall present a construction of proper biharmonic submanifolds in $\s^{n}$.
Let $M$ be a  submanifold of $\s^{n-1}(\frac{1}{\sqrt{2}})$.
Then $M$ can be seen as a submanifold of $\s^{n}$ and we have

\begin{theorem}[\cite{RCSMCO2},\cite{ELCO2}]\label{teo-min-bih}
Assume that $M$ is a submanifold of $\s^{n-1}(\frac{1}{\sqrt{2}})$. Then
$M$ is a proper biharmonic submanifold of  $\s^{n}$ if and only if it
is minimal in $\s^{n-1}(\frac{1}{\sqrt{2}})$.
\end{theorem}

Theorem~\ref{teo-min-bih} is a useful tool to construct examples of proper
biharmonic submanifolds. For instance, using a well known result
of H.B.~Lawson \cite{HBL},
we have

\begin{theorem}[\cite{RCSMCO2}]
There exist closed orientable embedded proper biharmonic
surfaces of arbitrary genus in  $\s^4$.
\end{theorem}

This shows the existence of an abundance of proper biharmonic surfaces in
$\s^4$, in contrast with the case of $\s^3$.

The biharmonic submanifolds that we have produced so far are all
pseudo-umbilical, i.e. $A=|H|^2 I$. We now want to give examples of biharmonic
submanifolds of $\s^n$ that are not of this type.

\noindent With this aim, let $n_1$, $n_2$ be two positive integers
such that $n=n_1+n_2$, and let $r_1$, $r_2$ be two positive real
numbers such that $r^2_1+r^2_2=1$. Let $M_1$ be a minimal
submanifold of $\s^{n_1}(r_1)$, of dimension $m_1$, with
$0<m_1<n_1$, and let $M_2$ be a minimal submanifold of
$\s^{n_2}(r_2)$, of dimension $m_2$, with $0<m_2<n_2$. We have:

\begin{theorem}[\cite{RCSMCO2}]
The manifold $M_1\times M_2$ is a proper biharmonic submanifold
of $\s^{n+1}$ if and only if $r_1=r_2=\frac{1}{\sqrt{2}}$ and
$m_1\neq m_2$.
\end{theorem}

If $M$ is a submanifold of $\s^n$ with $|H|=\cst$, then it is possible
to give a partial classification. In fact we have

\begin{theorem}[\cite{CO}]
Let $M$ be a submanifold of $\s^n$ such that $\vert
H\vert$ is constant.
\begin{itemize}
\item[a)] If $\vert H\vert>1$, then $M$ is never biharmonic.
\item[b)] If $\vert H\vert=1$, then $M$ is biharmonic if and only
if it is pseudo-umbilical and $\nabla^{\perp}H=0$, i.e. $M$ is a
minimal submanifold of  $\s^{n-1}(\frac{1}{\sqrt{2}})\subset\s^n$.
\end{itemize}
\end{theorem}

\noindent As an immediate consequence we have

\begin{corollary}[\cite{CO}]
If $M$ is a compact orientable hypersurface of $\s^n$ with $\vert
H\vert=1$, then $M$ is proper biharmonic if and only if
$M=\s^{n-1}(\frac{1} {\sqrt{2}})$.
\end{corollary}

We end this section presenting two classes of proper biharmonic curves of $\s^n$

\begin{proposition}[\cite{RCSMCO2}]\quad\\
\vspace*{-5mm}

\begin{enumerate}
\item[a)] The circles
$
\gamma(t)=\cos(\sqrt{2}t)c_1+\sin(\sqrt{2}t)c_2+c_4,
$
where $c_1$, $c_2$, $c_4$ are constant orthogonal vectors of $\r^{n+1}$ with
$\vert c_1\vert^2=\vert c_2\vert^2=\vert c_4\vert^2=\frac{1}{2}$,
are proper biharmonic curves of $k_1=1$.
\item[b)] The curves
$
\gamma(t)=\cos(at)c_1+\sin(at)c_2+\cos(bt)c_3+\sin(bt)c_4,
$
where $c_1$, $c_2$, $c_3$, $c_4$ are constant orthogonal vectors of $\r^{n+1}$ with $\vert
c_1\vert^2=\vert c_2\vert^2=\vert c_3\vert^2= \vert
c_4\vert^2=\frac{1}{2}$, and $a^2+b^2=2$, $a^2\neq b^2$, are
proper biharmonic of  $k_1^2=1-a^2b^2\in (0,1)$.
\end{enumerate}
\end{proposition}

\subsection{Biharmonic submanifolds in Sasakian space forms}

A ``generalization'' of Riemannian manifolds with constant sectional curvature is
that of  Sasakian space forms. First, recall that $(N,\eta,\xi,\varphi,g)$ is a {\it contact Riemannian manifold} if: $N$ is
a $(2r+1)-$dimensional manifold; $\eta$ is an one-form
satisfying $(d\eta)^r\wedge\eta\neq 0$; $\xi$ is the
vector field defined by
$\eta(\xi)=1$ and $d\eta(\xi,\cdot)=0$; $\varphi$ is an
endomorphism field; $g$ is a Riemannian metric on $N$ such that,
$\forall X,Y\in C(TN)$,
\begin{itemize}
\item $\varphi^2=-I+\eta\otimes \xi$
\item $g(\varphi X, \varphi Y)=g(X,Y)-\eta(X) \eta(Y), \quad g(\xi,\cdot)=\eta$
\item $d\eta(X,Y)=2g(X,\varphi Y)$.
\end{itemize}

A contact Riemannian manifold $(N,\eta,\xi,\varphi,g)$ is a {\it Sasaki manifold} if
$$
\left(\nabla_{X}\varphi\right)(Y)=g(X,Y)\xi-\eta(Y)X.
$$

If the sectional curvature is constant on all $\varphi$-invariant tangent
$2$-planes of $N$, then $N$ is called of
{\it constant holomorphic sectional curvature}.
Moreover, if a Sasaki manifold $N$ is connected, complete and 
of constant holomorphic sectional curvature,
then it is called a {\it Sasakian space form}. We have the following classification.

\begin{theorem}[\cite{JBFTLV}]
A simply connected three-dimensional Sasakian space form is isomorphic to one of
the following:
\begin{itemize}
\item[a)] the special unitary group $SU(2)$
\item[b)] the Heisenberg group $\h_3$
\item[c)] the universal covering group of $SL_2(\r)$.
\end{itemize}
\end{theorem}
In particular, a  simply connected three-dimensional Sasakian space form
of  constant holomorphic sectional curvature $1$ is isometric to $\s^3$.

In \cite{JI}, J. Inoguchi classified proper biharmonic Legendre
curves and Hopf cylinders in three-dimensional Sasakian space forms.
To state Inoguchi results we recall that:
\begin{itemize}
\item a curve $\gamma:I\to N$ parametrized by arc length
is {\it Legendre} if $\eta(\gamma')=0$;
\item a Hopf cylinder is $S_{\overline{\gamma}}=\pi^{-1}(\overline{\gamma})$,
where
$\pi:N \to \overline{N}=N/G$ is the projection of $N$ onto the orbit
space $\overline{N}$ determined by the action of the one-parameter group of
isometries generated by $\xi$, when the action is simply transitive.
\end{itemize}

\begin{theorem}
[\cite{JI}] Let $N^3(\epsilon)$ be a Sasakian space form of constant
holomorphic sectional curvature $\epsilon$ and $\gamma:I\to N$ a biharmonic
Legendre curve parametrized by arclength.
\begin{itemize}
\item[a)] If $\epsilon\leq 1$, then $\gamma$ is a Legendre geodesic.
\item[b)] If $\epsilon>1$, then $\gamma$ is a Legendre geodesic or a Legendre
helix of curvature $\sqrt{\epsilon-1}$.
\end{itemize}
\end{theorem}

\begin{theorem}
[\cite{JI}] Let $S_{\overline{\gamma}}\subset N^3(\epsilon)$ be a
biharmonic Hopf cylinder in a Sasakian space form.
\begin{itemize}
\item[a)] If $\epsilon\leq 1$, then $\overline{\gamma}$ is a geodesic.
\item[b)] If $\epsilon>1$, then $\overline{\gamma}$ is a geodesic or a
Riemannian circle of curvature $\overline{k}=\sqrt{\epsilon-1}$.
\end{itemize}
\noindent In particular, there exist proper biharmonic Hopf
cylinders in Sasakian space forms of holomorphic sectional curvature
greater than $1$.
\end{theorem}

T. Sasahara classified, in \cite{TS1}, proper biharmonic Legendre surfaces in
Sasakian space forms and, in the case when the ambient space is the unit 
$5-$dimensional sphere
$\s^5$, he obtained their explicit representations.

\begin{theorem}
[\cite{TS1}] Let $\phi:M^2\to \s^5$ be a proper biharmonic Legendre
immersion. Then the position vector field $x_0=x_0(u,v)$ of $M$ in
$\r^6$ is given by:
\begin{align*}
x_0(u,v)=\tfrac{1}{\sqrt{2}}&\big( \cos u,\sin u \sin(\sqrt{2}v),
-\sin u \cos(\sqrt{2}v),\\
& \sin u,\cos u \sin(\sqrt{2}v), -\cos u
\cos(\sqrt{2}v)\big).
\end{align*}
\end{theorem}

Other results on biharmonic Legendre curves and biharmonic anti-invariant 
surfaces in Sasakian space forms and $(k,\mu)$-manifolds were obtained in 
\cite{KARECMTS1,KARECMTS2}.

\section{Biharmonic Riemannian submersions}

In this section we discuss some examples of proper  biharmonic Riemannian
submersions. From the expression of the bitension field
\eqref{eq:bihar-submersion} we have immediately the following

\begin{theorem}[\cite{CO1}]\label{teo-sub}
Let $\phi:M\to N$ be a Riemannian submersion with basic, non-zero, tension field.
Then $\phi$ is proper biharmonic if:
\begin{itemize}
\item[a)] $^N\nabla\tau(\phi)=0$;
\item[b)] $\tau(\phi)$ is a unit Killing vector field on $N$.
\end{itemize}
\end{theorem}

Theorem~\ref{teo-sub} was used in \cite{CO1} to construct examples of proper
biharmonic Riemannian submersions. These examples are projections
$\pi:TM\to M$ from the tangent bundle of a Riemannian manifold
endowed with a ``Sasaki type'' metric.
Indeed, let $(M,g)$ be an $m$-dimensional Riemannian manifold and let
$\pi:TM\to M$ be its tangent bundle. We denote by $V(TM)$ the vertical
distribution on $TM$ defined by $V_v(TM)=
\ker d\pi_v$, $v\in TM$. We consider a nonlinear connection on
$TM$ defined by the distribution $H(TM)$ on $TM$, complementary to
$V(TM)$, i.e. $H_v(TM)\oplus V_v(TM)=T_v(TM)$, $v\in TM$.  For any induced local chart
$(\pi^{-1}(U);x^i,y^j)$ on $TM$ we have a local adapted frame in $H(TM)$ defined
by the local vector fields
$$
\frac{\delta}{\delta x^i}=\frac{\partial}{\partial x^i}-
N_i^j(x,y)\frac{\partial}{\partial y^j}, \quad i=1,\ldots,m,
$$
where the local functions $N^i_j(x,y)$ are the connection
coefficients of the nonlinear connection defined by $H(TM)$.
 If we endow $TM$ with the   Riemannian metric $S$  defined by
$$
S(X^V,Y^V)=S(X^H,Y^H)=g(X,Y), \quad
S(X^V,Y^H)=0,
$$
then the canonical projection $\pi:(TM,S)\to (M,g)$ is a Riemannian submersion.
(For more details on the metrics on the tangent bundle see, for example,
\cite{VO})
\noindent The biharmonicity  of the map $\pi$ depends on the choice of
the connection coefficients $N_i^j$. For suitable choices we have:

\begin{proposition}[\cite{CO1}]\quad\\
\vspace*{-5mm}
\begin{itemize}
\item[a)] Let $\xi$ be an unit Killing vector field and let
$N^i_j=(\Gamma^i_{jk}+\delta^i_j\xi_k+\delta^i_k\xi_j)y^k$ be
a projective change of the Levi-Civita connection $\nabla$
on $(M,g)$. Then
$\pi$ is a proper biharmonic map.
\item[b)] Let $\rho\in C^{\infty}(M),\,\rho\neq\cst$,  be an affine function
and let $N^i_j=(\Gamma^i_{jk}+\delta^i_j\alpha_k+\delta^i_k\alpha_j-
g_{jk}\alpha^i)y^k,\, \alpha_k=\frac{\partial\rho}{\partial x^k}$,
be a conformal change of the connection $\nabla$. Then $\pi$ is a proper biharmonic map.
\end{itemize}
\end{proposition}

\section{Biharmonic maps between Euclidean spaces}
Let $\phi:\r^m\to\r^n$, $\phi(x)=(\phi^1(x),\ldots,\phi^n(x)),\,x\in\r^m$ be a
smooth map. Then, the bitension field assumes the simple expression
$\tau_2(\phi)=(\Delta^2\phi^1,\ldots,\Delta^2\phi^n)$.
Thus, a map $\phi:\r^m\to\r^n$ is biharmonic if and only if its components
functions are biharmonic.

If we want proper solutions defined everywhere, then we can take polynomial
solutions of degree three. If we look for maps which are not defined
everywhere, then there are interesting classes of examples.
One of this can be described as follows.

A smooth map $\phi:\r^m\setminus\{0\}\to \r^m\setminus\{0\}$ is {\it axially symmetric}
if there exist a
map $\f:\s^{m-1}\to\s^{n-1}$ and a function
$\rho:(0,\infty)\to(0,\infty)$ such that, for $y\in
\r^m\setminus\{0\}$,
$$
\phi(y)=\rho(|y|)\f\Big(\frac{y}{|y|}\Big).
$$

Assume that the map $\varphi$ is not constant.
An axially symmetric map $\phi=\rho\times\f:
\r^m\setminus\{0\}\to\r^n\setminus\{0\}$ is harmonic if and
only if $\f$ is an eigenmap of eigenvalue $2k>0$ (see \cite{JEAR} for the definition of eigenmaps) and
\begin{equation}\label{eq:rho_armonic_axsymm1}
\rho(t)=c_1 t^{A_1}+c_2 t^{A_2} ,
\end{equation}
where $2A_{1,2}=-(m-2)\pm\sqrt{(m-2)^2+8k}$ and 
$c_1, c_2\geq 0$ with $c_1^2+c_2^2\neq 0$.

The biharmonicity of axially symmetric maps
 $\phi=\rho\times\f:
\r^m\setminus\{0\}\to\r^n\setminus\{0\}$ was discussed in \cite{ABSMCO},
where the authors give the following classification.

\begin{theorem}[\cite{ABSMCO}]
Let $\phi=\rho\times\f:
\r^m\setminus\{0\}\to\r^n\setminus\{0\}$ be an axially symmetric map
and assume that $\f$ is an eigenmap of eigenvalue $2k>0$.
\begin{enumerate}
\item[a)] If $\rho'=0$, then
\begin{itemize}
\item for $m\geq 4$, $\phi$ can not be biharmonic.
\item for $m=3$, $\phi$ is proper biharmonic if
and only if $\f$ is an eigenmap of homogeneous degree
$h=1$.
\item for $m=2$, $\phi$ is  proper biharmonic if
and only if $\f$ is an eigenmap of homogeneous degree
$h=2$.
\end{itemize}
\item[b)] If  $\rho'\neq 0$, then
 $\phi$ is proper biharmonic if and
only if
\begin{equation}\label{exp: ax_symm_r_m_to_r_n}
\rho(t)=\begin{cases}
    c_1t^3\ln t+c_2t\ln t+c_3\ln t+c_4, \quad {\rm when}\quad
    m=2\, {\rm and}\; k=\frac{1}{2}
    \\ \mbox{} \\
   \frac{c_1}{2(m+2A_1)}t^{A_1+2}+\frac{c_2}{2(m+2A_2)}t^{A_2+2}
   +c_3t^{A_1}+c_4t^{A_2} ,
\quad {\rm otherwise}.
  \end{cases}
\end{equation}
where 
$c_1^2+c_2^2\neq 0$ and $c_1, c_2, c_3, c_4\geq 0$.
\end{enumerate}
\end{theorem}

\begin{example}
An important class of axially symmetric diffeomorphisms of
$\r^m\setminus\{0\}$ is given by
$$
\phi:\r^m\setminus\{0\}\to \r^m\setminus\{0\},\quad
\phi(y)= y/|y|^\ell\,,\quad \ell\neq 0,1,
$$
which, for $\ell=2$, provides the well
known Kelvin transformation. For these maps,
$\rho(t)=1/t^{\ell-1}$ and $\f:\s^{m-1}\to\s^{m-1}$ is the
identity map. An easy computation shows that $\phi$ is
harmonic if and only if $m=\ell$.

Using \eqref{exp: ax_symm_r_m_to_r_n} it follows that
$\phi$ is proper biharmonic if and only if
$m=\ell+2$. For $\ell=2$ this result was first obtained in
\cite{PBDK}.

We also note that the proper biharmonic map $
\phi:\r^m\setminus\{0\}\to \r^m\setminus\{0\}$, $
\phi(y)=y/|y|^{m-2}$, is harmonic with respect to the conformal
metric on the domain given by
$\widetilde{g}=|y|^{\frac{4}{3-m}}g_{\textrm{can}}$. This property
is similar to that of the Kelvin transformation proved by B.~Fuglede
in \cite{BF}.
\end{example}

\section{Biharmonic maps and conformal changes}

\subsection{Conformal change on the domain}

Let $\phi:(M^m,g)\to (N^n,h)$ be a harmonic map.
Consider a conformal change of the domain metric, i.e. $\tilde{g}=e^{2\rho}g$
for some smooth function $\rho$.

\noindent If $m=2$, from the conformal invariance of the energy, the map
${\phi}:(M,\tilde{g})\to (N,h)$ remains harmonic.
If $m\neq 2$, then ${\phi}$ does not remain, necessarily, harmonic.
Therefore, it is reasonable to seek under what conditions on the
function $\rho$ the map ${\phi}:(M,\tilde{g})\to (N,h)$ is biharmonic.

This problem was attacked in \cite{PBDK}, where P.~Baird and D.~Kamissoko 
first proved the following general result.
\begin{proposition}[\cite{PBDK}]\label{pro-conf-cod}
Let $\phi:(M^m,g)\to (N^n,h),\,m\neq 2$, be a harmonic map. Let $\tilde{g}=e^{2\rho}g$
be a metric conformally equivalent to $g$. Then ${\phi}:(M,\tilde{g})\to (N,h)$
is biharmonic if and only if
\begin{eqnarray*}
&-\Delta d\phi(\grad \rho)+(m-6) \nabla_{\grad\rho}d\phi(\grad\rho)
+2(\Delta\rho-(m-4)|d\rho|^2)d\phi(\grad\rho)\\
&+\trace R^{N}(d\phi(\grad\rho),d\phi)d\phi=0\,.
\end{eqnarray*}
\end{proposition}

If $\phi:(M,g)\to (M,g)$ is the identity map
${\bf 1}$, we call a conformally equivalent metric
$\tilde{g}=e^{2\rho}g$, for which ${\bf 1}$ becomes biharmonic, a
{\it biharmonic metric} with respect to $g$.

\noindent Applying the maximum principle we have
\begin{theorem}
[\cite{PBDK}] Let $(M^m,g)$, $m\neq 2$, be a compact manifold of
negative Ricci curvature. Then there is no biharmonic metric
conformally related to $g$ other than a constant multiple of $g$.
\end{theorem}

There is a surprising connection between biharmonic metrics and isoparametric functions.
We recall
that a smooth function $f:M\to\r$ is called isoparametric if,
for each $x\in M$ where $\grad f_x\neq 0$, there are real
functions $\lambda$ and $\sigma$ such that
$$
|df|^2  = \lambda\circ f,\quad
\Delta f  = \sigma \circ f\,,
$$
on some neighbourhood of $x$.
The above mentioned link is provided by the following
\begin{theorem}[\cite{PBDK}]\label{teo-conf-dom-isop}
 Let $(M^m,g)$, $m\neq 2$, be an Einstein manifold. Let
$\tilde{g}=e^{2\rho}g$ be a biharmonic metric conformally equivalent
to $g$. Then the function $\rho:M\to \r$ is isoparametric.

\noindent Conversely, let $f:M\to \r$ be an isoparametric function, then away
from critical points of $f$, there is a reparametrization
$\rho=\rho\circ f$ such that $\tilde{g}=e^{2\rho}g$ is a biharmonic metric.
\end{theorem}

\subsection{Conformal change on the codomain}
Let $\phi:(M^m,g)\to (N^n,h)$ be a harmonic map.
Consider the ``dual problem'', i.e. a conformal change $\tilde{h}=e^{2\rho}h$ of the codomain metric.  In this case the analogous of
Proposition~\ref{pro-conf-cod} is more complicated and we shall  review only on
some special situations.

If ${\bf 1}:(M,g)\to(M,g)$ is the identity map, then it is proved, in  \cite{AB1},
that ${\bf 1}:(M,g)\to(M,e^{2\rho}g)$ is biharmonic if and only if
\begin{align*}
\trace\nabla^2\grad\rho+&(2\Delta\rho+(2-m)|\grad\rho|^2)\grad\rho
+\tfrac{6-m}{2}\grad(|\grad\rho|^2)\\
&+\ricci(\grad\rho)=0.
\end{align*}
This equation was used in \cite{AB1} to prove  similar results to
Theorem~\ref{teo-conf-dom-isop}, for the conformal change of the metric on the codomain.

In a similar setting, in \cite{CO3,CO4}, C.~Oniciuc constructed new examples of
biharmonic maps deforming the metric on a sphere. More precisely,
let $\s^n\subset\r^{n+1}$ be the  $n-$dimensional sphere endowed with the conformal
modified metric $e^{2\rho}\langle,\rangle$, where $\langle,\rangle$ is the canonical
metric on $\s^n$ and $\rho(x)=x^{n+1}$. Let $\s^{n-1}=\{x\in\s^n:
x^{n+1}=0\}$ be the equatorial sphere of $\s^n$. Then the inclusion
${\bf i}:(\s^{n-1},\langle,\rangle)\to
(\s^n,e^{2\rho}\langle,\rangle)
$
is a proper biharmonic map.

\noindent This result was generalized in
\begin{theorem}[\cite{CO3,CO4}]
Let $M$ be a minimal submanifold of
$(\s^{n-1},\langle,\rangle)$. Then $M$ is a proper biharmonic submanifold of
$(\s^n,e^{2\rho}\langle,\rangle)$.
\end{theorem}

Observe that even a geodesic $\gamma:I\to(N,h)$ will not
remain harmonic after a conformal change of the metric on $(N,h)$, unless the
conformal factor is constant.
As to biharmonicity of $\gamma$ we have the following.
\begin{theorem}
[\cite{ELSM}] Let $(N^{n},h)$ be a Riemannian manifold. Fix a point
$p\in N$ and let $f=f(r)$ be a non-constant function, depending
only on the geodesic distance $r$ from $p$, which is a solution of
the following ODE:
$$
 f'''(r)+3f''(r)f'(r)+f'(r)^3=0.
$$
Then any geodesic $\gamma:I\to(N,h)$ such that
$p\in\gamma(I)$ becomes a proper biharmonic curve $\gamma:I\to
(N,e^{2f}h)$.
\end{theorem}

For example, take $(N,h)=({\r}^{2},g=dx^2+dy^2)$ and
$f(r)= \ln{(r^{2} +1)}$,  where $r=\sqrt{x^2+y^2}$ is the distance
from the origin. Then any straight line on the flat $\r^2$ turns to
a biharmonic curve on $(\r^2,\bar{g}=(r^{2} +1)^2(dx^2+dy^2))$, which
is the metric, in local isothermal coordinates, of the Enneper
minimal surface.

\section{Biharmonic morphisms}
In analogy with the case of harmonic morphisms (see \cite{PBJCW}) the definition
of {\it biharmonic morphisms} can be formulated as follows.
\begin{definition}
A map $\phi:(M,g)\to (N,h)$ is a biharmonic morphism
if for any biharmonic function $f:U\subset N\to \r$, its pull-back
by $\phi$, $f\circ\phi:\phi^{-1}(U)\subset M\to \r$, is a biharmonic function.
\end{definition}

In ~\cite{ELYLO} E.~Loubeau and Y.-L.~Ou gave the characterization of the
biharmonic morphisms showing that a map is a biharmonic morphism if
and only if it is a horizontally weakly conformal biharmonic map and
its dilation satisfies a certain technical condition.

\noindent A more direct characterization is

\begin{theorem}
[\cite{YLO1,ELYLO}] A map $\phi:(M,g)\to (N,h)$ is a biharmonic morphism
if and only if there exists a function $\lambda:M\to \r$ such that
$$
\Delta^2(f\circ\phi)=\lambda^4\Delta^2(f)\circ\phi,
$$
for all functions $f:U\subset N\to \r$.
\end{theorem}

If $M$ is compact, the notion of biharmonic morphisms becomes trivial,
in fact we have

\begin{theorem}
[\cite{ELYLO}] Let $\phi:(M,g)\to (N,h)$ be a non-constant map. If
$M$ is compact, then $\phi$ is a biharmonic morphism if and only if
it is a harmonic morphism of constant dilation, hence a homothetic
submersion with minimal fibers.
\end{theorem}

In \cite{YLO2}, Y.-L.~Ou, using the theory of $p-$harmonic morphisms,
proved the following properties.

\begin{theorem}
[\cite{YLO2}] The radial projection $\phi:\r^m\setminus\{0\}\to
\s^{m-1}$, $\phi(x)=\frac{x}{|x|}$, is a biharmonic morphism if and
only if $m=4$.
\end{theorem}

\begin{theorem}
[\cite{YLO2}] The projection $\phi:M\times_{\beta^2}N\to (N,h)$,
$\phi(x,y)=y$, of a warped product onto its second factor is a
biharmonic morphism if and only if $1/\beta^2$ is a harmonic
function on $M$.
\end{theorem}

In the case of polynomial biharmonic
morphisms between Euclidean spaces there is a full classification.

\begin{theorem}
[\cite{YLO2}] Let $\phi:\r^m\to \r^n$ be a polynomial biharmonic
morphism, i.e. a biharmonic morphism whose component functions are
polynomials, with $m>n\geq2$. Then $\phi$ is an orthogonal
projection followed by a homothety.
\end{theorem}

\section{The second variation of biharmonic maps}

The second variation formula for the bienergy functional $E_2$ was
obtained, in a general setting, by G.Y.~Jiang in ~\cite{GYJ2}.
For biharmonic maps in Euclidean spheres, the second variation
formula takes a simpler expression.
\begin{theorem}[\cite{CO2}]
Let $\phi:(M,g)\to \s^n$ be a biharmonic map. Then the
Hessian of the bienergy $E_2$ at $\phi$ is given by
$$H(E_2)_{\phi}(V,W)=\int_M\langle I^{\phi}(V),W\rangle v_g,$$
where
\begin{eqnarray*}\label{eq:J1}
I^{\phi}(V)&=&\Delta(\Delta
V)+\Delta\{\trace\langle V,d\phi\cdot\rangle d\phi\cdot- \vert d\phi\vert^2V\}
+2\langle d\tau(\phi),d\phi\rangle V+\vert\tau(\phi)\vert^2V\nonumber\\
&&-2\trace\langle V,d\tau(\phi)\cdot\rangle d\phi\cdot
-2\trace\langle \tau(\phi),dV\cdot\rangle d\phi\cdot\nonumber\\
&&-\langle \tau(\phi),V\rangle \tau(\phi)+\trace\langle d\phi\cdot,\Delta V\rangle d\phi\cdot
\nonumber\\
&&+\trace\langle d\phi\cdot,\trace\langle V,d\phi\cdot\rangle d\phi\cdot\rangle d\phi\cdot
-2\vert d\phi\vert^2\trace\langle d\phi\cdot,V\rangle d\phi\cdot\nonumber\\
&&+2\langle dV,d\phi\rangle \tau(\phi)-\vert d\phi\vert^2\Delta V+ \vert
d\phi\vert^4 V.
\end{eqnarray*}
\end{theorem}

Although the expression of the operator $I$ is rather complicated,
in some particular cases it becomes easy to study.

\noindent In the instance when  $\phi$ is the identity map of $\s^n$,
$I^{\bf 1}$ has the expression
$$
I^{\bf 1}(V)=\Delta(\Delta V)-2(n-1)\Delta V+(n-1)^2V,
$$
and we can immediately  deduce
\begin{theorem}
[\cite{CO2}] The identity map ${\bf 1}:\s^n\to \s^n$ is biharmonic
stable and
\begin{enumerate}
\item[a)] if $n=2$, then $\nul({\bf 1})=6$;
\item[b)] if $n>2$, then $\nul({\bf 1})=\frac{n(n+1)}{2}$.
\end{enumerate}
\end{theorem}
A large class of biharmonic maps for which it is possible to study the Hessian
is obtained using the following generalization of Theorem~\ref{teo-min-bih}.
\begin{theorem}[\cite{ELCO2}]\label{teo-comp-gen}
Let $M$ be an orientable compact manifold and
${\bf i}:\s^{n-1}(\frac{1}{\sqrt{2}})\to\s^{n}$ the canonical inclusion. If
$\psi:M\to \s^{n-1}(\frac{1}{\sqrt{2}})$ is a non-constant map, then
$\phi={\bf i}\circ\psi:M\to \s^{n}$ is proper biharmonic if
and only if $\psi$ is harmonic and $e(\psi)$ is constant.
\end{theorem}

\begin{remark}\label{re-unstable}
All the biharmonic maps constructed using  Theorem~\ref{teo-comp-gen} are unstable.
To see this, let $\phi_t:\s^{n-1}\to \s^{n}$, $\phi_t(x)=(tx,\sqrt{1-t^2})$, $t\in
[0,1]$, the map defined in Section~\ref{section-sn}. Then
$$
\big(E_2(\phi_t)\big)''_{{t=\frac{1}{\sqrt{2}}}}
=-2(n-1)^2\vol(\s^{n-1})<0.
$$
\end{remark}

Thus the problem is to describe qualitatively their $\Index$ and $\nul$.

\noindent When $\psi$ is the identity map of $\s^{n-1}(\frac{1}{\sqrt{2}})$ we have

\begin{theorem}[\cite{ELCO1},\cite{ABCO}]\label{teo-ind-inc} The biharmonic $\Index$ of
the canonical inclusion ${\bf
i}:\s^{n-1}(\frac{1}{\sqrt 2})\to \s^{n}$ is exactly $1$, and its
$\nul$ is $\frac{n(n-1)}{2}+n$. 
\end{theorem}

\noindent When $\psi$ is the minimal generalized Veronese map we get

\begin{theorem}[\cite{ELCO1}]\label{teo-ind-ver} The biharmonic map derived
from the generalized
Veronese map $\psi:\s^m(\sqrt{\frac{m+1}{m}})\to
\s^{m+p}(\frac{1}{\sqrt{2}})$, $p=\frac{(m-1)(m+2)}{2}$, has $\Index$
at least $m+2$, when $m\leq 4$, and at least $2m+3$, when $m>4$.
\end{theorem}

In Theorem~\ref{teo-ind-inc} and~\ref{teo-ind-ver} the map $\psi$ was a
minimal immersion. We
shall consider now the case of harmonic Riemannian submersions, and
choose for $\psi$ the Hopf map.

\begin{theorem}
[\cite{ELCO2}] The $\Index$ of the biharmonic map
$\phi={\bf i}\circ\psi:\s^3(\sqrt{2})\to \s^3$ is at least $11$, while its $\nul$ is
bounded from below by $8$. 
\end{theorem}

\noindent We note that, for the above results, the authors described explicitly the
spaces where $I^{\phi}$ is negative definite or vanishes.

For the case of surfaces in Sasakian space forms, T. Sasahara,
considering a variational vector field parallel to $H$, gave a sufficient
condition for proper biharmonic Legendre submanifolds into an
arbitrary  Sasakian space form to be unstable.
This condition is expressed in terms of  the 
 mean curvature vector field and
of the second fundamental form of the submanifold. In particular

\begin{theorem}[\cite{TS3}] The biharmonic Legendre curves and surfaces in Sasakian
space forms are unstable.
\end{theorem}

\end{document}